\documentclass[11pt]{article}
\usepackage{mathrsfs}
\usepackage{amsmath}
\usepackage{amsmath,amssymb,amsfonts,amsthm,fancyhdr}
\usepackage{epsfig,graphicx,picins,picinpar,subfigure}
\usepackage{pstricks}
\usepackage{fancyvrb}
\usepackage[numbers,sort&compress]{natbib}
\begin{document}
\title{\textbf{Generalizations of an Ancient Greek Inequality about the Sequence of Primes}}
\author{\emph{\textbf{Shaohua Zhang}}$^{1,2}$}
\date{{\small 1 School of Mathematics, Shandong University,
Jinan,  Shandong, 250100, China \\
2 The key lab of cryptography technology and information security,
Ministry of Education, Shandong University, Jinan, Shandong, 250100,
China\\E-mail address: shaohuazhang@mail.sdu.edu.cn}}
 \maketitle

\begin{abstract}
In this note, we generalize an ancient Greek inequality about the
sequence of primes to the cases of arithmetic progressions even
multivariable polynomials with integral coefficients. We also refine
Bouniakowsky's conjecture [16] and Conjecture 2 in [22]. Moreover,
we give two remarks on conjectures in [22].

\vspace{3mm}\noindent \textbf{Keywords:} inequality, primes,
Euclid's second theorem, Dirichlet's theorem, Bouniakowsky's
conjecture

\vspace{3mm}\noindent \textbf{2000 MR  Subject Classification:}\quad
11Y41; 11A99

\end{abstract}

\section{INTRODUCTION}
In his \emph{Elements}, Euclid proved that prime numbers are more
than any assigned multitude of prime numbers. In other words, there
are infinitely many primes. For the details of proof, see [1,
Proposition 20,  Book 9]. Hardy and Wright [2] called this classical
result Euclid's second theorem. Hardy likes particularly Euclid's
proof. He [3] called it is "as fresh and significant as when it was
discovered---two thousand years have not written a wrinkle on it".
According to Hardy [3], "Euclid's theorem which states that the
number of primes is infinite is vital for the whole structure of
arithmetic. The primes are the raw material out of which we have to
build arithmetic, and Euclid's theorem assures us that we have
plenty of material for the task". Andr\'{e} Weil [4] also called
"the proof for the existence of infinitely many primes represents
undoubtedly a major advance......". Many people like Euclid's second
theorem. In his magnum opus \emph{History of the Theory of Numbers},
Dickson [5] gave the historical list of proofs of Euclid's second
theorem from Euclid (300 B.C.) to M\'{e}trod (1917). Ribenboim [6]
cited nine and a half proofs of Euclid's second theorem. The author
[7] cited fifteen new proofs.

\vspace{3mm} Based on Euclid's idea, people in Ancient Greek could
prove that for $n>1$, $\prod_{i=1}^{i=n}p_i>p_{n+1}$ since
$p_{n+1}\leq \prod_{i=1}^{i=n}p_i-1$, where $p_i$ represents the
$i^{th}$ prime. We call the inequality
$\prod_{i=1}^{i=n}p_i>p_{n+1}$ Ancient Greek inequality. In 1907,
Bonse [8] refined this inequality and proved that for $n\geq 4$,
$\prod_{i=1}^{i=n}p_i>p_{n+1}^2$ and for $n\geq 5$,
$\prod_{i=1}^{i=n}p_i>p_{n+1}^3$. This kind of inequalities has been
improved since then [9, 10]. Why are people interested in the
inequality between $\prod_{i=1}^{i=n}p_i$ and $p_{n+1}$? The main
reason is of that this kind of inequalities are closely related to
the famous Chebychev's function $\theta(x)=\sum_{p\leq x}\log p$.
And $\theta(x)\sim x \Longleftrightarrow \pi (x)\sim \frac{x}{\log
x}$ (The Prime Number Theorem).

\vspace{3mm}In a somewhat different direction, the aim of this note
is to generalize the ancient Greek inequality to the cases of
arithmetic progressions even multivariable polynomials with integral
coefficients. We noticed that $p_i$ can be viewed as the $i^{th}$
prime value of polynomial $f(x)=x$. Let $a$ and $b$ be integers with
$a\neq 0$, $b>0$ and $\gcd(a,b)=1$. Dirichlet's classical and most
important theorem states that  $f(x)=a+bx$ can represent infinitely
many primes. Denote the $i^{th}$ prime of the form $f(x)$ by
$P_{f,i}$. Naturally,  we want to prove that for every sufficiently
large integer $n$, $\prod_{i=1}^{i=n}P_{f,i}>P_{f,n+1}$, where
$f=a+bx$ with $a\neq 0$, $b>0$ and $\gcd(a,b)=1$. More generally, we
hope that if $f$ is a multivariable polynomial with integral
coefficients and $f$ can take infinitely many prime values, then
there is a constant $C$ such that when $n>C$,
$\prod_{i=1}^{i=n}P_{f,i}>P_{f,n+1}$. Thus, one could refine
Bouniakowsky's conjecture and so on. For the details, see Section 2.

\section{SOME THEOREMS AND CONJECTURES}
In this note, we always restrict that a $k$-variables polynomial
with integral coefficients is a map from $N^k$ to $Z$, where $k\in N
$ and $N$ is the set of all positive integers, $Z$ is the set of all
integers.

\vspace{3mm}Now, let's begin with Bertrand's and related problems in
arithmetic progressions. In 1845, Bertrand [5] verified for numbers
$<6000000$ that for any integer $n>6$ there exists at least one
prime between $n-2$ and $\frac{n}{2}$. In 1850, Chebychev [5] proved
that there exists a prime between $x$ and $2x-2$ for $x>3$.  In the
case of arithmetic progressions, Breusch [11], Ricci [12] and
Erd\"{o}s [13] proved respectively that for $n \geq 6$, positive
integer, there is always a prime $p$ of the form $6n + 1$, and one
of the form $6n -1$, such that $n < p < 2n$. This implies
immediately that the following Theorem 1 and Theorem 2.

\vspace{3mm}\noindent {\bf Theorem 1:~~}%
Let $f(x)=6x+1$. Then $P_{f,1}=7$,  $P_{f,2}=13$, $P_{f,3}=19$,....
And for $n>1$, $\prod_{i=1}^{i=n}P_{f,i}>P_{f,n+1}$.

\vspace{3mm}\noindent {\bf Theorem 2:~~}%
Let $f(x)=6x-1$. Then $P_{f,1}=5$,  $P_{f,2}=11$, $P_{f,3}=17$,....
And for $n>1$, $\prod_{i=1}^{i=n}P_{f,i}>P_{f,n+1}$.

\vspace{3mm} In 1941, Molsen [14] proved  (1) for $n\geq 199$, the
interval $n<p\leq{\textstyle\frac 8{7}}n$ always contains a prime of
each of the forms $3x+1$, $3x-1$; (2) for $n\geq 118$, the interval
$n<p\leq{\textstyle\frac 4{3}}n$ always contains a prime of each of
the forms $12x+1,12x-1,12x+5,12x-5$. Based on Molsen's work, it is
not difficult to  prove that the following theorems.

\vspace{3mm}\noindent {\bf Theorem 3:~~}%
Let $f(x)=3x+1$. Then $P_{f,1}=7$,  $P_{f,2}=13$, $P_{f,3}=19$,....
And for $n>1$, $\prod_{i=1}^{i=n}P_{f,i}>P_{f,n+1}$.

\vspace{3mm}\noindent {\bf Theorem 4:~~}%
Let $f(x)=3x-1$. Then $P_{f,1}=2$,  $P_{f,2}=5$, $P_{f,3}=11$,....
And for $n>2$, $\prod_{i=1}^{i=n}P_{f,i}>P_{f,n+1}$.

\vspace{3mm}\noindent {\bf Theorem 5:~~}%
Let $f(x)=4x+1$. Then $P_{f,1}=5$,  $P_{f,2}=13$, $P_{f,3}=17$,....
And for $n>1$, $\prod_{i=1}^{i=n}P_{f,i}>P_{f,n+1}$.

\vspace{3mm}\noindent {\bf Theorem 6:~~}%
Let $f(x)=4x-1$. Then $P_{f,1}=3$,  $P_{f,2}=7$, $P_{f,3}=11$,....
And for $n>1$, $\prod_{i=1}^{i=n}P_{f,i}>P_{f,n+1}$.

\vspace{3mm} Let $a$ and $b$ be integers with $a\neq 0$, $b>0$ and
$\gcd(a,b)=1$. In 1896, Ch. de la Vall\'{e}e-Poussin [15] proved
that $\sum _{p\equiv a(\mod b),p\leq x}\log p$ equals
$\frac{x}{\varphi (b)}$ asymptotically. Therefore, for  every
sufficiently large integer $n$, $\sum _{i=1}^{i=n+1}\log P_{f,i}$
equals $\frac{P_{f,n+1}}{\varphi (b)}$ asymptotically. Clearly,
$\frac{P_{f,n+1}}{\varphi (b)}>2\log P_{f,n+1}$. It shows
immediately that the following Theorem 7 holds.

\vspace{3mm}\noindent {\bf Theorem 7:~~}%
 Let $a$ and $b$ be integers with $a\neq 0$, $b>0$ and
$\gcd(a,b)=1$. And let $f(x)=a+bx$. Then there is a constant $C$
depending on $a$ and $b$ such that when $n>C$,
$\prod_{i=1}^{i=n}P_{f,i}>P_{f,n+1}$.

\vspace{3mm}Based on the aforementioned theorems, also based on
Bateman-Horn's heuristic asymptotic formula [17], we give a
strengthened form of Bouniakowsky's conjecture [16] which can be
viewed as a refinement of special form of Schinzel-Sierpinski's
Conjecture [18] as follows:

\vspace{3mm}\noindent{\bf  Conjecture 1:~~}%
If $f(x)$ is an irreducible polynomial with integral coefficients,
positive leading coefficient, and there does not exist any integer
$n>1$ dividing all the values $f(k)$ for every integer $k$, then
$f(x)$ represents  primes for infinitely many $x$, moreover, there
is a constant $C$ such that when $n>C$,
$\prod_{i=1}^{i=n}P_{f,i}>P_{f,n+1}$.

\vspace{3mm}  Conjecture 1 can be deduced by Bateman-Horn's formula.
Next, we will try to generalize Conjecture 1 to the cases of
multivariable polynomials with integral coefficients. Firstly, we
have the following theorems:

\vspace{3mm}\noindent {\bf Theorem 8 [19]:~~}%
Let $f(x,y)=x^2+y^2+1$. Then there is a constant $C$ such that when
$n>C$, $\prod_{i=1}^{i=n}P_{f,i}>P_{f,n+1}$.

\vspace{3mm}\noindent {\bf Theorem 9 [20]:~~}%
Let $f(x,y)=x^2+y^4$. Then there is a constant $C$ such that when
$n>C$, $\prod_{i=1}^{i=n}P_{f,i}>P_{f,n+1}$.

\vspace{3mm}\noindent {\bf Theorem 10 [21]:~~}%
Let $f(x,y)=x^3+2y^3$. Then there is a constant $C$ such that when
$n>C$, $\prod_{i=1}^{i=n}P_{f,i}>P_{f,n+1}$.

\vspace{3mm}By the aforementioned idea and theorems, one could
strengthen a special form of Conjecture 2 in [22] as follows:

\vspace{3mm}\noindent{\bf Conjecture 2:~~}%
Let $f(x_1,...,x_k)$ be a multivariable polynomial with  integral
coefficients, if there is a positive integer $c$ such that for every
positive integer $m\geq c$, there exists an integral point
$(y_1,...,y_k)$ such that $f(y_1,...,y_k)>1$ is in $Z_m^*=\{x\in
N|\gcd (x,m)=1, x\leq m\}$, and there exists an integral point
$(z_1,...,z_k)$ such that $f(z_1,...,z_k)\geq c$ is prime, then
$f(x_1,...,x_k)$ represents  primes for infinitely many integral
points $(x_1,...,x_k)$. Moreover, there is a constant $C$ such that
when $n>C$, $\prod _{i=1}^{i=n}P_{f,i}>P_{f,n+1}$.

\vspace{3mm}\noindent {\bf Remark 1:~~}%
Conjecture 2 implies that a special case of Conjecture 1 in [22].
Namely, if $f(x_1,...,x_k)$ is a multivariable polynomial with
integral coefficients, and represents primes for infinitely many
integral points $(x_1,...,x_k)$, then there is always a constant $c$
such that for every positive integer $m>c$, there exists an integral
point $(y_1,...,y_k)$ such that $f(y_1,...,y_k)>1$ is in $Z_m^*$. In
fact, by Conjecture 2, we know that there is a constant $C$ such
that when $n>C$, $\prod_{i=1}^{i=n}P_{f,i}>P_{f,n+1}$. Let $C<k\leq
C+1$ and let $c=P_{f,k}$.  When $m>c$,  we can assume that $c\leq
P_{f,k+h}\leq m<P_{f,k+h+1}$ with $h\geq 0$. If for some $1\leq
r\leq k+h$, $\gcd (P_{f,r},m)=1$, then there exists an integral
point $(y_1,...,y_k)$ such that $f(y_1,...,y_k)=P_{f,r}$ is in
$Z_m^*$. If for any $1\leq r\leq k+h$, $\gcd (P_{f,r},m)>1$, then
$m\geq \prod_{i=1}^{i=k+h}P_{f,i}>P_{f,k+h+1}$ since $C<k\leq k+h$
and $\gcd (P_{f,i},P_{f,j})=1$ for $i\neq j$. It is a contradiction.

\vspace{3mm}\noindent {\bf Remark 2:~~}%
Conjecture 1 can not be extended to arbitrary number-theoretic
functions without a proviso. For example, let
$$h(n)=\left\{
\begin{array}{c}
p_1=2, n=1 \\
p_2=3, n=2\\
......\\
\mbox{the least prime of the form } k\times \prod_{i=1}^{i=n-1}p_i+1, n\geq 2\\
\end{array}
\right. .$$

\vspace{3mm} Clearly, for any positive integer $n$,
$\prod_{i=1}^{i=n}P_{h,i}<P_{h,n+1}$, where $P_{h,i}=h(i)$ is the
$i^{th}$ prime value of the function $h(n)$.

\vspace{3mm} By this example, one also can find that Conjecture 1 in
[22] can not be extended to arbitrary number-theoretic functions
without a proviso. In fact, if there is such a constant $c$, then
there is always a positive integer $k$  such that $c<h(k)$. Let
$m=\prod_{i=1}^{i=k}P_{h,i}$. Clearly, in this case, $c<m\leq
h(k+1)-1<h(k+1)$ and there does not exist any positive integer $y$
such that $h(y)>1$ is in $Z_m^* $. Otherwise, $y\geq k+1$. It is
impossible since $m<h(k+1)$.

\vspace{3mm}Let  $s\geq1$ and $k\geq 1$ be integers. Let
$f_1(x_1,...,x_k),...,f_s(x_1,...,x_k)$ be multivariable polynomials
with  integral coefficients. We also assume  that $f_1(x_1,...,x_k),
..., f_s(x_1,...,x_k)$ represent simultaneously primes for
infinitely many integral points $(x_1,...,x_k)$. Denote the set of
integral points $(x_1,...,x_k)$ such that
$f_1(x_1,...,x_k),...,f_s(x_1,...,x_k)$ are primes by $X$. Let
$\beta_{f,1}=\prod_{i=1}^{i=s}f_i(X_1)$, where $X_1\in X$ such that
the norm $||(f_1(X_1),...,f_s(X_1))||$ is the least. Let
$\beta_{f,2}=\prod_{i=1}^{i=s}f_i(X_2)$, where $X_2\in X$ such that
$\gcd(\beta_{f,1},\beta_{f,2})=1$,  $||(f_1(X_1),...,f_s(X_1))||<
||(f_1(X_2),...,f_s(X_2))||\leq ||(f_1(X_0),...,f_s(X_0))||$ with
$X_0\in X,X_0\neq X_1, X_0\neq X_2$ and
$\gcd(\prod_{i=1}^{i=s}f_i(X_0), \beta_{f,1})=1$, ... Let
$\beta_{f,j}=\prod_{i=1}^{i=s}f_i(X_j)$, where $X_j\in X$ such that
for any $1\leq r\leq j-1$,
$\gcd(\beta_{f,j},\prod_{i=1}^{i=s}f_i(X_1)\times
...\times\prod_{i=1}^{i=s}f_i(X_{j-1}))=1$,
$||(f_1(X_r),...,f_s(X_r))||< ||(f_1(X_j),...,f_s(X_j))||\leq
||(f_1(X_0),...,f_s(X_0))||$ with $X_0\in X, X_0\neq X_1, X_0\neq
X_2,...,X_0\neq X_j$ and
$\gcd(\prod_{i=1}^{i=s}f_i(X_0),\prod_{i=1}^{i=s}f_i(X_1)\times
...\times\prod_{i=1}^{i=s}f_i(X_{j-1}))=1$, ... Clearly,
$\gcd(\beta_{f,i},\beta_{f,j})=1$ for $i\neq j$. Notice that
pairwise distinct primes are pairwise relatively prime. The sequence
of primes $\{p_i\}$ has a beautiful property: if any integral
sequence $1<a_1<...a_n<...$ with $\gcd(a_i,a_j)=1$ for $i\neq j$,
then $p_i\leq a_i$ for any positive integer $i$. For the proof of
this property, see Appendix. Therefore, like the $i^{th}$ prime
$p_i$, $\beta_{f,i}$ can be viewed as the $i^{th}$ "desired prime
number". Thus, one could give a strengthened form of Conjecture 2 in
[22] as follows:

\vspace{3mm}\noindent{\bf  Conjecture 3:~~}%
Let $f_1(x_1,...,x_k)$, ..., $f_s(x_1,...,x_k)$ be multivariable
polynomials with  integral coefficients, if there is a positive
integer $c$ such that for every positive integer $m\geq c$, there
exists an integral point $(y_1,...,y_k)$ such that
$f_1(y_1,...,y_k)>1,...,f_s(y_1,...,y_k)>1$ are all in $Z_m^*
=\{x\in N|\gcd (x,m)=1, x\leq m\}$, and there exists an integral
point $(z_1,...,z_k)$ such that $f_1(z_1,...,z_k)\geq
c,...,f_s(z_1,...,z_k)\geq c$ are all primes, then
$f_1(x_1,...,x_k)$, ..., $f_s(x_1,...,x_k)$  represent
simultaneously primes for infinitely many integral points
$(x_1,...,x_k)$. Moreover, there is a constant $C$ such that when
$n>C$, $\prod_{i=1}^{i=n}\beta_{f,i}>\beta_{f,n+1}$.

\section{CONCLUSIONS}
In this note, we generalized an ancient Greek inequality about the
sequence of primes to the cases of arithmetic progressions. By
Bateman-Horn's heuristic asymptotic formula and also based on the
work of Motohashi Yoichi, Friedlander John, Iwaniec Henryk,
Heath-Brown, and so on,  we refined Bouniakowsky's conjecture and
Conjecture 2 in [22]. Knuth called Euclid's Algorithm the granddaddy
of all algorithms. Based on the work in this note, one can see that
the Ancient Greek inequality about the sequence of primes also is
the granddaddy of the inequalities about the sequence of some kind
special kinds of primes.

\section{ACKNOWLEDGEMENTS}
Thank my advisor Professor Xiaoyun Wang for her valuable help. Thank
for the Institute Advanced Study in Tsinghua University for
providing me with excellent conditions. This work was partially
supported by the National Basic Research Program (973) of China (No.
2007CB807902) and the Natural Science Foundation of Shandong
Province (No. Y2008G23).

\section{REFERENCES}

\vspace{3mm}\noindent[1] Thomas Little Heath, The Thirteen Books of
the Elements, translated from the text of Heiberg with introduction
and commentary, Cambridge Univ. Press, Cambridge (1926).

\vspace{3mm}\noindent[2] Hardy G. H. and Wright E. M., An
Introduction to the Theory of Numbers, Oxford: The Clarendon Press,
(1938).

\vspace{3mm}\noindent[3] G. H. Hardy, A Mathematician's Apology,
Cambridge University Press, (1940).

\vspace{3mm}\noindent[4] Andr\'{e} Weil, Number theory, an approach
through history from Hammurapi to Legendre. Birkh\"{a}ser, Boston,
Inc., Cambridge, Mass., (1984).

\vspace{3mm}\noindent[5] Leonard Eugene Dickson, History of the
Theory of Numbers, Volume I: Divisibility and Primality, Chelsea
Publishing Company, New York, (1952).

\vspace{3mm}\noindent[6] Paulo Ribenboim, The book of prime number
records, Springer-Verlag, New York, (1988).

\vspace{3mm}\noindent[7] Shaohua Zhang, Euclid's Number-Theoretical
Work, available at: \\http://arxiv.org/abs/0902.2465

\vspace{3mm}\noindent[8] H. Bonse, \"{U}er eine bekannte Eigenschaft
der Zahl 30 und ihre Verallgemeinerung, Arch. Math. Phys., 12,
292-295, (1907).

\vspace{3mm}\noindent[9] P\'{o}sa Lajos, \"{U}ber eine Eigenschaft
der Primzahlen, (Hungarian) Mat. Lapok, 11, 124-129, (1960).

\vspace{3mm}\noindent[10] Panaitopol Lauren\c{t}iu, An inequality
involving prime numbers, Univ. Beograd. Publ. Elektrotehn. Fak. Ser.
Mat., 11, 33-35, (2000).

\vspace{3mm}\noindent[11] R. Breusch, Zur Verallgemeinerung der
Bertrandschen Postulates dass zwischen x und 2x stets Primzahlen
liegen, Math. Z., 34, 505-526, (1932).

\vspace{3mm}\noindent[12] G. Ricci, Sul teorema di Dirichlet
relativo alla progresione aritmetica, Boll. Un. Mat. Ital., 12,
304-309,  (1933).

\vspace{3mm}\noindent[13] P. Erd\"{o}s, \"{U}ber die Primzahlen
gewisser arithmetischen Reihen, Math. Zeit., 39, 473-491, (1935).

\vspace{3mm}\noindent[14] K. Molsen, Zur Verallgemeinerung der
Bertrandschen Postulates. Deutsche Math. 6, 248-256, (1941).

\vspace{3mm}\noindent[15] Ch. de la Vall\'{e}e-Poussin, Recherches
analytiques sur la th\'{e}orie des nombres (3 parts), Ann. Sec. Sci.
Bruxelles, 20, 183-256; 361-397, (1896).

\vspace{3mm}\noindent[16] Bouniakowsky, V.,  Nouveaux
th\'{e}or\`{e}mes relatifs \`{a} la distinction des nombres premiers
et \`{a}la d\'{e} composition des entiers en facteurs, Sc. Math.
Phys., 6, 305-329, (1857).

\vspace{3mm}\noindent[17] Paul T. Bateman and Roger A. Horn, A
heuristic asymptotic formula concerning the distribution of prime
numbers, Math. Comp., Vol. 16, No.79, 363-367, (1962).

\vspace{3mm}\noindent[18] A. Schinzel and W. Sierpinski, Sur
certaines hypotheses concernant les nombres premiers, Acta Arith., 4
(1958), 185-208, Erratum 5, 259, (1958).

\vspace{3mm}\noindent[19] Motohashi Yoichi, On the distribution of
prime numbers which are of the form $x^2+y^2+1$, Acta Arith., 16,
351-363, (1969/1970).

\vspace{3mm}\noindent[20] Friedlander John, Iwaniec Henryk, The
polynomial $x^2+y^4$ captures its primes, Ann. of Math., (2) 148,
no. 3, 945-1040, (1998).

\vspace{3mm}\noindent[21] Heath-Brown D. R. Primes represented by
$x^3+2y^3$, Acta Math., 186, No. 1, 1-84, (2001).

\vspace{3mm}\noindent[22] Shaohua Zhang, On the Infinitude of Some
Special Kinds of Primes, available at:
http://arxiv.org/abs/0905.1655

\section{APPENDIX}
In this appendix, we prove the following theorem 11:

\vspace{3mm}\noindent{\bf  Lemma 1:~~}%
$\pi (n)$ is the largest among the cardinality of all sub-sets in
which each element exceeds 1 and pairwise distinct elements are
pairwise relatively prime of $\{1,2,...,n\}$, where $\pi (x)$
represents the number of primes less than or equal to $x$.

\vspace{3mm}\noindent{\bf  Proof:~~}%
Easy. Let $S$ be a sub-set of $\{1,2,...,n\}$ such that in $S$, each
element exceeds 1 and pairwise distinct elements are pairwise
relatively prime. Denote the cardinality of $S$ by $|S|$. If
$|S|>\pi (n)$, then $\prod_{x\in S}x $ has at least $\pi (n)+1$
distinct prime divisors. This implies that there must be an element
$a\in S$ such that $a\geq p_{\pi (n)+1}>n$. It is a contradiction
since $S\subseteq \{1,2,...,n\}$. This complets the proof of Lemma
1.

\vspace{3mm}\noindent{\bf  Theorem 11:~~}%
If any integral sequence $1<a_1<...a_n<...$ with $\gcd(a_r,a_j)=1$
for $r\neq j$, then $p_i\leq a_i$ for any positive integer $i$,
where $p_i$ is the $i^{th}$ prime.

\vspace{3mm}\noindent{\bf  Proof:~~}%
Easy. For any positive integer $i$, we consider the set
$S=\{a_1,...,a_i\}$. By known condition, we have $\gcd(a_r,a_j)=1$
for $r\neq j$. Namely, in $S$, each element exceeds 1 and pairwise
distinct elements are pairwise relatively prime. So, by Lemma 1,
$i\leq \pi (a_i)$. It shows that $p_i\leq a_i$. Therefore, Theorem
11 holds. This completes the proof of Theorem 11.

\clearpage
\end{document}